\documentclass[11pt,reqno]{amsart}

\usepackage[a4paper,left=0.75in,right=0.75in,top=1in,bottom=1in]{geometry}
\usepackage[T1]{fontenc}
\usepackage{lmodern}
\usepackage{amsmath,amssymb,amsthm,mathtools}
\usepackage{microtype}
\usepackage[colorlinks=true,linkcolor=blue,citecolor=blue,urlcolor=blue,bookmarks=false]{hyperref}

\numberwithin{equation}{section}
\theoremstyle{plain}
\newtheorem{theorem}{Theorem}[section]
\newtheorem{lemma}[theorem]{Lemma}
\newtheorem{proposition}[theorem]{Proposition}
\newtheorem{corollary}[theorem]{Corollary}
\theoremstyle{remark}

\newcommand{\dd}{\,\mathrm d}
\newcommand{\RR}{\mathbb R}
\newcommand{\Q}{\mathcal Q}
\newcommand{\D}{{D_\theta}}
\newcommand{\ub}{u_b}
\newcommand{\Rr}{\mathcal R}
\newcommand{\Ss}{\mathcal S}

\begin{document}

\title[Robin monotonicity along regular polygons]{Monotonicity of the first Robin eigenvalue along regular polygons of fixed area or perimeter}
\author{Yijing Chen}
\address{School of Mathematics and Statistics, HNP-LAMA, Central South University, Changsha, Hunan 410083, P. R. China}
\email{chenyijingpde@qq.com}
\author{Zhuo Cheng}
\address{School of Mathematics and Statistics, HNP-LAMA, Central South University, Changsha, Hunan 410083, P. R. China}
\email{zhuo@outlook.cz}
\date{}

\begin{abstract}
We establish strict monotonicity of the first Robin eigenvalue along planar regular polygons of fixed area or fixed perimeter, with the Robin parameter $\alpha$ in the boundary condition $\partial_\nu u+\alpha u=0$ held fixed. Under either normalization, the eigenvalue strictly decreases with the number of sides for $\alpha>0$, strictly increases for $\alpha<0$, and is identically zero for $\alpha=0$.

The proof is based on a continuous family of mixed Robin--Neumann problems on fundamental right triangles. For negative parameters, a stress-tensor Rellich identity gives the sign of the angular derivative under both normalizations. For positive parameters in the fixed-area setting, a polar Rellich decomposition, directional monotonicity, and a transverse power-concavity property of the ground state yield a uniform estimate that absorbs a sign-changing boundary term into the bulk contribution for polygons with at least four sides. The triangle-to-square comparison follows from explicit secular equations. The positive-parameter result for fixed perimeter is then deduced by scaling. The proof is entirely analytical and applies to every $N\ge3$ and every real Robin parameter. Passing to the Robin-to-Dirichlet limit and using real analyticity of the limiting mixed Dirichlet--Neumann branch, we also obtain strict monotonicity of the first Dirichlet eigenvalue under both geometric normalizations.
\end{abstract}

\maketitle

\section{Introduction}
\label{sec:intro}

\subsection{Background and motivation}
For a bounded connected Lipschitz domain $\Omega\subset\RR^2$ and a real parameter $\alpha\in\RR$, the principal eigenvalue of the Robin Laplacian is characterized variationally by
\begin{equation}\label{eq:robin-rayleigh}
 \lambda_1^\alpha(\Omega)=
 \min_{v\in H^1(\Omega)\setminus\{0\}}
 \frac{\displaystyle\int_\Omega|\nabla v|^2+
       \alpha\int_{\partial\Omega}v^2}
      {\displaystyle\int_\Omega v^2},
\end{equation}
corresponding to the boundary value problem $-\Delta u=\lambda u$ in $\Omega$ with $\partial_\nu u+\alpha u=0$ on $\partial\Omega$, where $\nu$ denotes the outward unit normal. As $\alpha$ varies over the real line, \eqref{eq:robin-rayleigh} interpolates between the Neumann problem ($\alpha=0$, where $\lambda_1^0(\Omega)=0$), the elastically supported membrane ($\alpha>0$), and the Dirichlet problem in the limit $\alpha\to+\infty$, while extending for $\alpha<0$ to attractive boundary interactions arising in reaction--diffusion models and surface superconductivity; see the survey \cite{BFK}.

When $\Omega$ ranges over all planar domains of fixed area $A>0$, the extremal behavior of $\lambda_1^\alpha(\Omega)$ dates back to Rayleigh's conjecture for the clamped membrane, proved independently by Faber and Krahn via symmetrization \cite{Henrot}. Extending the Faber--Krahn inequality to $\alpha>0$ required methods beyond Schwarz rearrangement, because the level sets of a Robin eigenfunction meet $\partial\Omega$ transversely. Bossel \cite{Bossel} in dimension two and Daners \cite{Daners} in arbitrary dimension proved that the disk $B_A$ of area $A$ uniquely minimizes $\lambda_1^\alpha(\Omega)$ for every $\alpha>0$; an alternative free-discontinuity proof for general Lipschitz domains was later given by Bucur and Giacomini \cite{BucurGiacomini}. When $\alpha<0$, the negative boundary energy favors larger perimeter, and Bareket \cite{Bareket} conjectured that $B_A$ maximizes $\lambda_1^\alpha(\Omega)$ among domains of area $A$. Freitas and Krej\v{c}i\v{r}\'{\i}k \cite{FreitasKrejcirik} showed that this fails for large $|\alpha|$, as thin annuli eventually surpass the disk, although $B_A$ remains the maximizer in two dimensions for $|\alpha|$ sufficiently small.

A classical refinement of isoperimetric inequalities, initiated by P\'olya and Szeg\H{o} in \cite{PolyaSzego}, asks how geometric functionals behave on polygons and, in particular, along sequences of regular $N$-gons with fixed area or fixed perimeter as the number of sides $N\ge3$ increases toward the disk. Write $P_{N,A}$ for the regular $N$-gon of area $A>0$ and $Q_{N,L}$ for the regular $N$-gon of perimeter $L>0$. P\'olya and Szeg\H{o} proved via Steiner symmetrization that $P_{3,A}$ and $P_{4,A}$ minimize the first Dirichlet eigenvalue $\lambda_1^D$ among all triangles and quadrilaterals of area $A$ (the corresponding polygonal Faber--Krahn problem remains open for $N\ge5$), and proposed that $\lambda_1^D(P_{N,A})$ should decrease strictly with $N$. Even for quantities closely tied to conformal geometry, such discrete monotonicity is delicate: the analogous conjecture of P\'olya and Szeg\H{o} for logarithmic capacity was settled only half a century later by Solynin and Zalgaller \cite{SolyninZalgaller}. For eigenvalues of the Laplacian, distinct regular polygons of equal area satisfy no inclusion relation, Steiner and Schwarz symmetrizations destroy regularity, and for $N\ge5$ the eigenfunctions are no longer trigonometric combinations.

The two natural geometric normalizations lead to related but distinct spectral problems. Along the fixed-area sequence $\{P_{N,A}\}_{N\ge3}$, perimeter decreases strictly with $N$, whereas along the fixed-perimeter sequence $\{Q_{N,L}\}_{N\ge3}$, area increases strictly with $N$. Both observations suggest that the principal Robin eigenvalue should decrease for positive parameters and increase for negative parameters, in agreement with the sign of the boundary term in \eqref{eq:robin-rayleigh}. However, since the Robin scaling law is
\[
 \lambda_1^\alpha(r\Omega)=r^{-2}\lambda_1^{r\alpha}(\Omega),\qquad r>0,
\]
the two normalizations are not equivalent by homogeneity alone when the physical coupling parameter $\alpha$ is held fixed. In particular, for negative parameters the gradient and boundary terms respond to dilation with competing signs.

\subsection{Main results}
We prove strict monotonicity of the first Robin eigenvalue in the number of sides for every real parameter $\alpha\neq0$, under either fixed area or fixed perimeter. In both cases the physical coupling parameter $\alpha$ remains fixed. The principal result is the fixed-area theorem:

\begin{theorem}\label{thm:main}
Let $A>0$, and let $P_{N,A}\subset\RR^2$ be the regular $N$-gon of area $A$. For every integer $N\ge3$ and every $\alpha\in\RR$,
\begin{equation}\label{eq:main-inequalities}
\begin{cases}
\lambda_1^\alpha(P_{N+1,A})<\lambda_1^\alpha(P_{N,A}),&\alpha>0,\\[1.5mm]
\lambda_1^\alpha(P_{N+1,A})=\lambda_1^\alpha(P_{N,A})=0,&\alpha=0,\\[1.5mm]
\lambda_1^\alpha(P_{N+1,A})>\lambda_1^\alpha(P_{N,A}),&\alpha<0.
\end{cases}
\end{equation}
\end{theorem}

The proof is entirely analytical and uniform in $N\ge3$ and $\alpha\in\RR$: it relies neither on asymptotic expansions for large $N$ or extreme values of $\alpha$, nor on numerical computation. Passing to the limit $\alpha\to+\infty$ recovers $\lambda_1^D(P_{N,A})$, and although pointwise limits preserve only weak inequalities, the real analyticity of the underlying continuous-aperture Dirichlet branch upgrades weak monotonicity to strict inequality. We thus obtain strict Dirichlet monotonicity as a direct consequence of the positive Robin result.

\begin{corollary}\label{cor:dirichlet}
Let $A>0$, and let $B_A\subset\RR^2$ denote the disk of area $A$. For every integer $N\ge3$,
\begin{equation}\label{eq:dirichlet-monotonicity}
 \lambda_1^D(P_{N+1,A})<\lambda_1^D(P_{N,A}),
 \qquad
 \lim_{N\to\infty}\lambda_1^D(P_{N,A})=\lambda_1^D(B_A).
\end{equation}
\end{corollary}

The same monotonicity laws hold at fixed perimeter. For positive Robin parameters, this follows from Theorem~\ref{thm:main} and the strict decrease of the eigenvalue under dilation. For negative parameters, the stress-tensor identity used in the proof of Theorem~\ref{thm:main} also controls the fixed-perimeter branch. The Dirichlet case follows from Corollary~\ref{cor:dirichlet} and homogeneity.

\begin{corollary}[Fixed-perimeter monotonicity]\label{cor:perimeter}
Let $L>0$, and let $Q_{N,L}\subset\RR^2$ be the regular $N$-gon of perimeter $L$. For every integer $N\ge3$ and every $\alpha\in\RR$,
\begin{equation}\label{eq:perimeter-inequalities}
\begin{cases}
\lambda_1^\alpha(Q_{N+1,L})<\lambda_1^\alpha(Q_{N,L}),&\alpha>0,\\[1.5mm]
\lambda_1^\alpha(Q_{N+1,L})=\lambda_1^\alpha(Q_{N,L})=0,&\alpha=0,\\[1.5mm]
\lambda_1^\alpha(Q_{N+1,L})>\lambda_1^\alpha(Q_{N,L}),&\alpha<0.
\end{cases}
\end{equation}
Moreover, if $B_{L^2/(4\pi)}$ denotes the disk of perimeter $L$, then
\begin{equation}\label{eq:perimeter-dirichlet}
 \lambda_1^D(Q_{N+1,L})<\lambda_1^D(Q_{N,L}),
 \qquad
 \lim_{N\to\infty}\lambda_1^D(Q_{N,L})=\lambda_1^D(B_{L^2/(4\pi)}).
\end{equation}
\end{corollary}

\subsection{Outline of the proof}
Because $P_{N,A}$ is tiled by $2N$ congruent right triangles of apex angle $\pi/N$ and area $A/(2N)$, dihedral symmetry and simplicity of the principal eigenvalue identify $\lambda_1^\alpha(P_{N,A})$ with the value at $\theta=\pi/N$ of the first eigenvalue $\Lambda(\theta)$ of a mixed boundary problem on the right triangle
\[
 D_h(\theta)=\{(s,t)\in\RR^2:0<s<h(\theta),\ 0<t<s\tan\theta\},
 \qquad h(\theta)=\left(\frac{A\theta}{\pi\tan\theta}\right)^{1/2},
\]
with $\partial_s u_h+\alpha u_h=0$ on the vertical leg $s=h(\theta)$ and homogeneous Neumann data on the two rays $t=0$ and $t=s\tan\theta$.

In Section~\ref{sec:branch}, we differentiate $\Lambda(\theta)$ for $\theta\in(0,\pi/3]$ at fixed $A$ and $\alpha$. Up to a positive trigonometric factor, $\Lambda'(\theta)$ has the sign of
\[
 J_h=(1-c)Y_h-(1+c)X_h+\frac{1-c}{2}\alpha B_h,
 \qquad c=\frac{\sin\theta\cos\theta}{\theta}\in(0,1),
\]
where $X_h=\int_{D_h}|(u_h)_t|^2$, $Y_h=\int_{D_h}|(u_h)_s|^2$, and $B_h=\int_0^{h\tan\theta}u_h(h,t)^2\dd t$ for the $L^2$-normalized positive eigenfunction $u_h$. After justifying corner regularity up to second derivatives via reflection and an exponential gauge at the Robin--Neumann vertex (Lemma~\ref{lem:corner-regularity}), we dispose of the case $\alpha<0$ immediately (Proposition~\ref{prop:negative}): integrating the divergence-free stress tensor of $u_h$ against $(h-s)e_s$ eliminates the Robin leg and expresses $Y_h+\frac{\alpha}{2}B_h$ as a strictly negative boundary integral on the sloping Neumann edge, giving $J_h<0$ throughout $(0,\pi/3]$.

The bulk of the paper treats $\alpha>0$. At any fixed $\theta$, rescaling $D_h(\theta)$ to the unit-height triangle $\D=\{0<s<1,\ 0<t<qs\}$ ($q=\tan\theta$) with dimensionless coupling $\beta=\alpha h(\theta)>0$ and eigenvalue $\mu=h(\theta)^2\Lambda(\theta)>0$ yields $h^2 J_h = J := (1-c)Y-(1+c)X+\frac{1-c}{2}\beta B$ for the rescaled ground state $u$. In Section~\ref{sec:polar}, we combine an angular integration by parts with a polar Rellich identity to obtain the exact decomposition (Proposition~\ref{prop:polar})
\[
 J = P_1 + P_2 + B_\beta,
 \qquad
 B_\beta = \beta\int_0^q \left[2ct\arctan t-\frac{1-c}{2}\right]\ub(t)^2\dd t,
\]
where $\ub(t)=u(1,t)$. Unlike the negative-parameter case, $B_\beta$ does not have a definite sign: its bracketed kernel changes sign from negative near $t=0$ to positive near $t=q$, while $\ub(t)$ turns out to be decreasing and hence weights the negative region more heavily. Thus $B_\beta$ must be absorbed into the bulk integral $P_1$.

Sections~\ref{sec:directions}--\ref{sec:close} carry out this absorption for $0<\theta\le\pi/4$ (which covers all $N\ge4$). In Section~\ref{sec:directions}, using the ground-state representation of the shifted quadratic form together with $\cos(2\theta)\ge0$, we establish the pointwise signs $u_s<0$, $u_t<0$, $s u_t - t u_s > 0$, $u_{st}>0$, and $u_{ss}-u_{tt}\ge0$ in $\D$. These inequalities ensure in particular that $P_1>0$ and $P_2>0$.

To control how sharply $\ub(t)$ can concentrate near $t=0$, Section~\ref{sec:concavity} adapts power-concavity maximum principles in the spirit of Brascamp--Lieb \cite{BrascampLieb} and Korevaar \cite{Korevaar} to a directional setting adapted to the wedge angle. We prove that $\partial_{tt}(u^a)\le0$ on $\overline\D$ for $a=2\cos^2\theta$ (Proposition~\ref{prop:power}). The choice $a=2\cos^2\theta$ cancels the cubic gradient term in the oblique derivative identity for $(\log u)_{tt}+a(\log u)_t^2$ along the sloping Neumann side, allowing the strong maximum principle and the Hopf lemma to rule out positive maxima on $\overline\D$. Concavity of $\ub^a$ on $[0,q]$ and tangential comparison at $t=q$ then imply that $\ub(t)^2/w(t)$ is nondecreasing for the explicit weight $w(t)=[1+2\beta q(q-t)/(1+q^2)]^{1+q^2}$ (Corollary~\ref{cor:weight}).

In Section~\ref{sec:close}, one-dimensional variation of constants along horizontal segments using $u_{ss}+(\mu/2)u\ge0$ and $-u_{tt}\ge(\mu/2)u$ bounds $P_1$ from below by $\mu^2 M_3\int_0^q f(t)\ub(t)^2\dd t$ with a nonnegative kernel $f(t)$, and comparison against a separable trial function yields $\mu^2 M_3 > \beta^2/4$ for all $\beta>0$ and $\theta\in(0,\pi/4]$. The combined boundary integrand $[K(t)+(\beta/4)f(t)]\ub(t)^2$ has a kernel with a single sign change from negative to positive on $(0,q)$. Monotonicity of $\ub^2/w$ transfers positivity down to the explicit weight $w$, and a second single-crossing comparison in the exponent $1+q^2\in(1,2]$ reduces the integral to an explicit positive cubic polynomial, proving $\Lambda'(\theta)>0$ on $(0,\pi/4]$.

Beyond $\theta=\pi/4$, the sign of $\cos(2\theta)$ changes, so the preceding argument does not cover the triangle-to-square step. In Section~\ref{sec:final}, we compare the equilateral-triangle secular equation \cite{RudnickWigman} (with a self-contained verification for the ground state) with that of the square, proving $\lambda_1^\alpha(P_{3,A})>\lambda_1^\alpha(P_{4,A})$ for all $\alpha>0$. This completes the proof of Theorem~\ref{thm:main}; the Robin-to-Dirichlet limit and analytic perturbation of the mixed Dirichlet--Neumann branch then give Corollary~\ref{cor:dirichlet}.

Finally, Section~\ref{sec:perimeter} establishes Corollary~\ref{cor:perimeter}. The fixed-perimeter positive-parameter comparison is a consequence of the fixed-area theorem and a dilation argument. For negative parameters, we change the height of the mixed-triangle branch to $h_L(\theta)=L\theta/(2\pi\tan\theta)$. Its derivative is again a positive factor times a linear combination of $Y_h+\alpha B_h/2$ and $-X_h$, with nonnegative coefficients, so the same Rellich identity gives strict monotonicity without any additional ground-state estimates.

\section{The mixed-triangle branch, regularity, and the case \texorpdfstring{$\alpha\le0$}{alpha <= 0}}
\label{sec:branch}

Fix $A>0$ and $\alpha\in\RR$. For $0<\theta\le\pi/3$, define
\[
 q=\tan\theta,\qquad
 h=h(\theta)=\left(\frac{A\theta}{\pi\tan\theta}\right)^{1/2},\qquad
 a=a(\theta)=h(\theta)\tan\theta,
\]
and consider the right triangle
\[
 D_h(\theta)=\{(s,t)\in\RR^2:0<s<h,\ 0<t<s\tan\theta\}
\]
of area $\frac12 ha = \frac{A\theta}{2\pi}$. We impose the Robin condition $\partial_s u+\alpha u=0$ on the vertical side $\Gamma_{R,h}=\{s=h,\ 0<t<a\}$ and homogeneous Neumann conditions on the horizontal and sloping sides. Let $\Lambda(\theta)$ denote the first eigenvalue of this mixed problem, and let $u_h$ be its unique positive eigenfunction normalized by $\int_{D_h(\theta)}u_h^2=1$. When $\theta=\pi/N$ ($N\ge3$), even reflection across the symmetry axes of $P_{N,A}$ and simplicity of $\lambda_1^\alpha(P_{N,A})$ give
\begin{equation}\label{eq:polygon-link}
 \Lambda(\pi/N)=\lambda_1^\alpha(P_{N,A}).
\end{equation}

Set
\[
 X_h=\int_{D_h(\theta)}|(u_h)_t|^2,\qquad
 Y_h=\int_{D_h(\theta)}|(u_h)_s|^2,\qquad
 B_h=\int_0^a u_h(h,t)^2\dd t,
\]
so that $\Lambda(\theta)=X_h+Y_h+\alpha B_h$. Pulling back $D_h(\theta)$ to the fixed reference triangle $D_*=\{(\xi,\eta):0<\eta<\xi<1\}$ via $(s,t)=(h\xi,a\eta)$ transforms the quadratic form and $L^2$-norm into
\[
 \mathfrak a_\theta[v]
 =\frac{a}{h}\int_{D_*}|v_\xi|^2
  +\frac{h}{a}\int_{D_*}|v_\eta|^2
  +\alpha a\int_0^1 v(1,\eta)^2\dd\eta,
 \qquad
 \mathfrak m_\theta[v]=ha\int_{D_*}v^2,
\]
defined on the $\theta$-independent domain $H^1(D_*)$. Since $\Lambda(\theta)$ is simple and the coefficients depend smoothly on $\theta$, the implicit function theorem applied to the normalized weak equation shows that $\theta\mapsto(\Lambda(\theta),u_h\circ(h\cdot,a\cdot))$ is $C^1$ on $(0,\pi/3]$ (with one-sided derivative at $\theta=\pi/3$). Differentiating $\mathfrak a_\theta[v_\theta]=\Lambda(\theta)\mathfrak m_\theta[v_\theta]$ at the normalized eigenfunction $v_\theta( \xi,\eta)=u_h(h\xi,a\eta)$ eliminates the eigenfunction derivative by the weak Euler--Lagrange equation and yields
\[
 \Lambda'(\theta)
 =\left(\frac{a'}{a}-\frac{h'}{h}\right)(Y_h-X_h)
  +\alpha\frac{a'}{a}B_h
  -\Lambda(\theta)\left(\frac{h'}{h}+\frac{a'}{a}\right).
\]
From $\log h=\frac12(\log\theta-\log\tan\theta+\text{const})$ and $a=h\tan\theta$, one computes
\[
 \frac{h'}{h}=\frac12\left(\frac1\theta-\frac1{\sin\theta\cos\theta}\right),
 \qquad
 \frac{a'}{a}=\frac12\left(\frac1\theta+\frac1{\sin\theta\cos\theta}\right).
\]
Substituting $\Lambda=X_h+Y_h+\alpha B_h$ and writing $c=c_\theta:=\frac{\sin\theta\cos\theta}{\theta}\in(0,1)$, we obtain
\begin{equation}\label{eq:branch-derivative}
 \Lambda'(\theta)=\frac{J_h}{\sin\theta\cos\theta},\qquad
 J_h=(1-c)Y_h-(1+c)X_h+\frac{1-c}{2}\alpha B_h.
\end{equation}

To analyze the sign of $J_h$ at a fixed angle $\theta\in(0,\pi/3]$, we rescale $D_h(\theta)$ to the unit-height triangle
\[
 \D=\{(s,t)\in\RR^2:0<s<1,\ 0<t<qs\},\qquad q=\tan\theta,
\]
with sides $\Gamma_R=\{s=1,\ 0<t<q\}$, $\Gamma_0=\{0<s<1,\ t=0\}$, and $\Gamma_\theta=\{0<s<1,\ t=qs\}$. Setting
\[
 \beta=\alpha h,\qquad \mu=h^2\Lambda(\theta),\qquad u(s,t)=h\,u_h(hs,ht),\qquad \ub(t)=u(1,t),
\]
the function $u>0$ satisfies $\int_\D u^2=1$ and solves
\begin{equation}\label{eq:unit-pde}
 -\Delta u=\mu u \quad\text{in }\D,\qquad
 u_s+\beta u=0 \quad\text{on }\Gamma_R,\qquad
 \partial_n u=0 \quad\text{on }\Gamma_0\cup\Gamma_\theta,
\end{equation}
where $n$ denotes the outward unit normal to $\partial\D$. Defining
\[
 X=\int_\D u_t^2=h^2X_h,\qquad
 Y=\int_\D u_s^2=h^2Y_h,\qquad
 B=\int_0^q \ub(t)^2\dd t=hB_h,
\]
we have
\begin{equation}\label{eq:J-unit}
 J := (1-c)Y-(1+c)X+\frac{1-c}{2}\beta B = h^2 J_h.
\end{equation}
We emphasize that this rescaling is performed solely at a frozen angle $\theta$ to determine the sign of $J_h$; when $\Lambda(\theta)$ is differentiated in \eqref{eq:branch-derivative}, it is $\alpha$, not $\beta$, that is held constant.

Before proceeding, we record the corner regularity of $u$ and the ground-state representation of the associated quadratic form, valid for every $\theta\in(0,\pi/3]$ and every $\beta\in\RR$.

\begin{lemma}[Corner regularity and positivity]\label{lem:corner-regularity}
For each fixed $\theta\in(0,\pi/3]$ and $\beta\in\RR$, the normalized positive eigenfunction $u$ of \eqref{eq:unit-pde} belongs to $H^{3+\varepsilon}(\D)\cap C^{2,\delta}(\overline\D)$ for some $\varepsilon,\delta>0$, and satisfies $\min_{\overline\D}u>0$. In particular, its first and second partial derivatives lie in $H^1(\D)\cap C(\overline\D)$, and the stress tensor
\[
 \mathbb T=\nabla u\otimes\nabla u-\frac12\bigl(|\nabla u|^2-\mu u^2\bigr)I
\]
belongs to $W^{1,1}(\D)$ with integrable boundary traces on $\partial\D$.
\end{lemma}

\begin{proof}
Interior and open-edge smoothness follow from standard elliptic regularity. At the right-angle corner $(1,0)$, even reflection across $\Gamma_0$ produces a domain with a smooth vertical Robin boundary, so $u$ is smooth up to $(1,0)$. At the origin $(0,0)$, both sides carry homogeneous Neumann data and meet at angle $\theta\le\pi/3$; the leading nonconstant Neumann angular mode behaves like $r^{\pi/\theta}\cos(\pi\varphi/\theta)$ with $\pi/\theta\ge3$, giving $H^{3+\varepsilon}\cap C^{2,\delta}$ regularity near $(0,0)$.

At the remaining vertex $(1,q)$, $\Gamma_R$ and $\Gamma_\theta$ meet at angle $\omega=\pi/2-\theta<\pi/2$. Choose a constant vector $b\in\RR^2$ such that $b\cdot n=\beta$ on $\Gamma_R$ and $b\cdot n=0$ on $\Gamma_\theta$. Then $z=e^{b\cdot x}u$ satisfies homogeneous Neumann boundary conditions on both edges of the corner and solves a second-order elliptic equation with constant coefficients. Even reflection across $\Gamma_\theta$ eliminates the odd mode of exponent $\pi/(2\omega)$, leaving $\pi/\omega>2$ as the lowest nonconstant angular exponent. Corner elliptic theory \cite{Grisvard} therefore yields $z\in H^{3+\varepsilon}\cap C^{2,\delta}$ near $(1,q)$ for small $\varepsilon,\delta>0$, and hence the same holds for $u$.

Strict positivity in $\D$ and on the open edges follows from the strong maximum principle and the Hopf boundary point lemma. At the vertices, strict positivity follows from the local Harnack inequality applied to $z$ (or from the corner expansion, since if $u$ vanished at a vertex, the first nonvanishing homogeneous Neumann mode would change sign inside the sector). Thus $\min_{\overline\D}u>0$, and $\mathbb T\in W^{1,1}(\D)$ follows immediately from $u\in H^2(\D)$.
\end{proof}

By Lemma~\ref{lem:corner-regularity}, differentiating the boundary conditions along the open edges and passing to the vertices by continuity gives the exact corner relations
\begin{equation}\label{eq:corner-values}
 \begin{aligned}
 &\nabla u(0,0)=0,\qquad D^2u(0,0)=-\frac\mu2 u(0,0)I,\\
 &u_s(1,q)=-\beta u(1,q),\qquad u_t(1,q)=-\beta q\,u(1,q),\qquad u_{st}(1,q)=\beta^2 q\,u(1,q).
 \end{aligned}
\end{equation}
Furthermore, the shifted quadratic form on $H^1(\D)$,
\begin{equation}\label{eq:Q-def}
 \Q[v]:=\int_\D|\nabla v|^2+\beta\int_0^q v(1,t)^2\dd t-\mu\int_\D v^2,
\end{equation}
admits the ground-state representation
\begin{equation}\label{eq:ground-transform}
 \Q[v]=\int_\D u^2\left|\nabla\left(\frac{v}{u}\right)\right|^2\ge0,
 \qquad
 \Q[v]=0\ \Longleftrightarrow\ v\in\operatorname{span}\{u\}.
\end{equation}
Indeed, since $\min_{\overline\D}u>0$, expanding $|\nabla(v/u)|^2$ and testing \eqref{eq:unit-pde} against $v^2/u\in H^1(\D)$ establishes \eqref{eq:ground-transform} for bounded $v\in H^1(\D)$, and general $v\in H^1(\D)$ follow by density.

We can already settle the nonpositive parameter regime across the full angular range.

\begin{proposition}[Monotonicity for $\alpha\le0$]\label{prop:negative}
If $\alpha<0$, then $\Lambda'(\theta)<0$ for every $\theta\in(0,\pi/3]$, and consequently $\lambda_1^\alpha(P_{N+1,A})>\lambda_1^\alpha(P_{N,A})$ for all $N\ge3$. If $\alpha=0$, then $\lambda_1^0(P_{N,A})=0$ for all $N\ge3$.
\end{proposition}

\begin{proof}
When $\alpha=0$, constant functions are admissible minimizers in \eqref{eq:robin-rayleigh}, giving $\lambda_1^0(P_{N,A})=0$. Assume now that $\alpha<0$, so $\beta=\alpha h<0$. Using the constant function $1$ as a trial function in the Rayleigh quotient for \eqref{eq:unit-pde} gives $\mu\le\beta q/\operatorname{Area}(\D)=2\beta<0$.

Consider the divergence-free tensor $\mathbb T$ from Lemma~\ref{lem:corner-regularity} and the smooth vector field $V=(1-s)e_s$ on $\D$. Since $\nabla V=-e_s\otimes e_s$ and $\mathbb T\in W^{1,1}(\D)$, the Gauss--Green theorem on the Lipschitz domain $\D$ applies directly without corner cutoffs:
\[
 \int_{\partial\D}(\mathbb T n)\cdot V\dd\ell
 =\int_\D \mathbb T:\nabla V
 =-\int_\D\left(\frac12 u_s^2-\frac12 u_t^2+\frac\mu2 u^2\right)
 =-\frac12(Y-X+\mu)=-\left(Y+\frac\beta2 B\right),
\]
where the last step uses $\mu=X+Y+\beta B$. On $\Gamma_R$, $V=0$; on $\Gamma_0$, $n=-e_t$ and $\partial_n u=0$, so $(\mathbb T n)\cdot e_s=0$. On $\Gamma_\theta$, $\partial_n u=0$ and the outward normal $n=(-\sin\theta,\cos\theta)$ has $e_s$-component $-\sin\theta$, whence $(\mathbb T n)\cdot V = \frac12\sin\theta\,(1-s)(|\nabla u|^2-\mu u^2)$. Therefore,
\[
 Y+\frac\beta2 B = -\frac{\sin\theta}{2}\int_{\Gamma_\theta}(1-s)\bigl(|\nabla u|^2-\mu u^2\bigr)\dd\ell < 0,
\]
because $\mu<0$ and $u>0$ on $\overline\D$. Since $0<c<1$ and $X\ge0$, \eqref{eq:J-unit} gives
\[
 J=(1-c)\left(Y+\frac\beta2 B\right)-(1+c)X<0,
\]
so $\Lambda'(\theta)<0$ on $(0,\pi/3]$ by \eqref{eq:branch-derivative}. Integrating $\Lambda'$ over $[\pi/(N+1),\pi/N]$ and using \eqref{eq:polygon-link} yields $\lambda_1^\alpha(P_{N+1,A})>\lambda_1^\alpha(P_{N,A})$ for every $N\ge3$.
\end{proof}

\section{Polar Rellich decomposition of the derivative form}
\label{sec:polar}

We now turn to the repulsive case $\alpha>0$ (so $\beta>0$ and $\mu>0$), working on the unit-height triangle $\D$. Introduce polar coordinates $(r,\varphi)$ centered at the apex $(0,0)$:
\[
 s=r\cos\varphi,\qquad t=r\sin\varphi,\qquad
 0<\varphi<\theta,\quad 0<r<\rho(\varphi):=\sec\varphi.
\]
Put $U(r,\varphi)=u(r\cos\varphi,r\sin\varphi)$ and define the radial and angular gradient components together with the auxiliary angular weight $g$:
\[
 \Rr=-U_r,\qquad \Ss=\frac{U_\varphi}{r},\qquad
 g(\varphi)=\sin\varphi\cos\varphi-c\varphi.
\]
Note that $g(0)=g(\theta)=0$, $g'(\varphi)=\cos(2\varphi)-c$, and concavity of $\sin(2\varphi)$ on $[0,\theta]$ implies $g(\varphi)>0$ for all $\varphi\in(0,\theta)$. In polar coordinates, the Cartesian derivatives read
\begin{equation}\label{eq:cartesian-polar}
 u_s=-\Rr\cos\varphi-\Ss\sin\varphi,\qquad
 u_t=-\Rr\sin\varphi+\Ss\cos\varphi.
\end{equation}

\begin{proposition}[Exact polar identity]\label{prop:polar}
For every $\theta\in(0,\pi/3]$ and $\beta\in\RR$, the form $J$ in \eqref{eq:J-unit} admits the exact representation
\begin{equation}\label{eq:J-decomp}
 J=P_1+P_2+B_\beta,
\end{equation}
where
\begin{align}
 P_1 &=-2\mu\int_0^\theta g(\varphi)\sec^4\varphi
       \int_0^1 s^2u(s,s\tan\varphi)u_t(s,s\tan\varphi)\dd s\dd\varphi, \label{eq:P1}\\
 P_2 &=2c\int_\D\Ss\bigl[(\tan\varphi+\varphi\sec^2\varphi)\Rr-\Ss\bigr]\dd A, \label{eq:P2}\\
 B_\beta &=\beta\int_0^q\left[2ct\arctan t-\frac{1-c}{2}\right]\ub(t)^2\dd t. \label{eq:Bbeta}
\end{align}
\end{proposition}

\begin{proof}
Using \eqref{eq:cartesian-polar} and $g'(\varphi)=\cos(2\varphi)-c$, we first split the bulk gradient difference in $J$ as
\[
 (1-c)Y-(1+c)X = I_1+I_2,
\]
where
\[
 I_1=\int_\D g'(\varphi)(\Rr^2-\Ss^2)\dd A,\qquad
 I_2=2\int_\D\Ss\bigl(\sin(2\varphi)\,\Rr-c\Ss\bigr)\dd A.
\]
Set $M=\int_\D g U U_\varphi\dd A$, $N=\int_\D g\tan\varphi\,U^2\dd A$, and define $E(\varphi)=\int_0^{\rho(\varphi)}(U_r^2-r^{-2}U_\varphi^2)r\dd r$. Differentiating $E(\varphi)$ using $\rho'(\varphi)=\rho(\varphi)\tan\varphi$, integrating the radial term $\int_0^\rho 2U_r U_{r\varphi}r\dd r$ by parts in $r$, and invoking the polar equation $r^{-1}(r U_r)_r+r^{-2}U_{\varphi\varphi}=-\mu U$ together with $\nabla u(0,0)=0$, we find
\[
 E'(\varphi)=\rho^2\tan\varphi\,(\Rr^2-\Ss^2)-2\rho^2\Rr\Ss+2\mu\int_0^\rho U U_\varphi r\dd r,
\]
where unintegrated functions of $r$ are evaluated on $\Gamma_R$ at $r=\rho(\varphi)$. Integrating $g'(\varphi)E(\varphi)$ by parts in $\varphi$ and using $g(0)=g(\theta)=0$ gives
\[
 I_1 = T_1 - 2\mu M,\qquad
 T_1 := \int_0^\theta g\rho^2\bigl[2\Rr\Ss-\tan\varphi\,(\Rr^2-\Ss^2)\bigr]\dd\varphi.
\]

Next, apply the Gauss--Green theorem on $\D$ to the divergence-free stress tensor $\mathbb T\in W^{1,1}(\D)$ contracted against the vector field $V=2g(\varphi)\tan\varphi\,r e_r$, which is smooth on $\overline\D\setminus\{(0,0)\}$ and Lipschitz at the origin. In the orthonormal polar frame $\{e_r,e_\varphi\}$, a direct computation yields
\[
 \mathbb T:\nabla V = 2\mu g\tan\varphi\,U^2 - 2(g\tan\varphi)'\Rr\Ss.
\]
The boundary flux $(\mathbb T n)\cdot V$ vanishes on $\Gamma_0\cup\Gamma_\theta$ because $V$ is radial while $\partial_n u=0$ and $n=\pm e_\varphi$ there. On $\Gamma_R$, the outward line element is $n\dd\ell=(\rho e_r-\rho' e_\varphi)\dd\varphi=\rho(e_r-\tan\varphi\,e_\varphi)\dd\varphi$, so the divergence theorem gives
\[
 T_2 := \int_0^\theta g\tan\varphi\,\rho^2\bigl[\Rr^2-\Ss^2+\mu U^2+2\tan\varphi\,\Rr\Ss\bigr]\dd\varphi
 = 2\mu N - 2\int_\D (g\tan\varphi)'\Rr\Ss\dd A.
\]
Along $\Gamma_R$, the Robin condition $u_s=-\beta u$ and \eqref{eq:cartesian-polar} read $\Rr\cos\varphi+\Ss\sin\varphi=\beta U$ and $-\Rr\sin\varphi+\Ss\cos\varphi=u_t$. Solving this orthogonal system for $\Rr$ and $\Ss$ on $\Gamma_R$ gives
\[
 2\Rr\Ss-2\tan\varphi\,(\Rr^2-\Ss^2)-2\tan^2\varphi\,\Rr\Ss = 2\beta U u_t\sec^2\varphi,
\]
so subtracting $T_2$ from $T_1$ yields
\[
 T_1-T_2 = \int_0^\theta g\rho^2\bigl[2\beta U u_t\sec^2\varphi-\mu\tan\varphi\,U^2\bigr]\dd\varphi.
\]
Meanwhile, $\sin(2\varphi)-(g\tan\varphi)' = c(\tan\varphi+\varphi\sec^2\varphi)$, so combining $I_1$ and $I_2$ produces
\[
 (1-c)Y-(1+c)X = 2\mu(N-M) + P_2 + 2\beta\int_0^\theta g\sec^4\varphi\,U u_t\dd\varphi - \mu\int_0^\theta g\rho^2\tan\varphi\,U^2\dd\varphi.
\]

To identify the $\mu$-terms with $P_1$, put $F(\varphi)=\int_0^{\rho(\varphi)}U(r,\varphi)^2 r\dd r$. Since $F'(\varphi)=\rho^2\tan\varphi\,U(\rho,\varphi)^2+2\int_0^\rho U U_\varphi r\dd r$, the two terms multiplying $\mu$ combine into $\mu\int_0^\theta g(\varphi)[2\tan\varphi\,F(\varphi)-F'(\varphi)]\dd\varphi$. Substituting $r=s\sec\varphi$ gives $F(\varphi)=\sec^2\varphi\int_0^1 s\,u(s,s\tan\varphi)^2\dd s$, whence
\[
 F'(\varphi)-2\tan\varphi\,F(\varphi) = 2\sec^4\varphi\int_0^1 s^2 u(s,s\tan\varphi)u_t(s,s\tan\varphi)\dd s,
\]
recovering $P_1$ in \eqref{eq:P1}. Finally, along $\Gamma_R$, $\frac{\dd}{\dd\varphi}\bigl[\ub(\tan\varphi)^2\bigr]=2\sec^2\varphi\,U(\rho(\varphi),\varphi)u_t(1,\tan\varphi)$. Integrating by parts in $\varphi$ and changing variables to $t=\tan\varphi$ gives
\[
 2\beta\int_0^\theta g\sec^4\varphi\,U u_t\dd\varphi
 =-\beta\int_0^\theta(g\sec^2\varphi)'\ub(\tan\varphi)^2\dd\varphi
 =-\beta\int_0^q \bigl[(1-c)-2ct\arctan t\bigr]\ub(t)^2\dd t,
\]
since $(g\sec^2\varphi)'=\sec^2\varphi[(1-c)-2c\varphi\tan\varphi]$. Adding $\frac{1-c}{2}\beta\int_0^q\ub(t)^2\dd t$ yields \eqref{eq:Bbeta} and completes the proof.
\end{proof}

Proposition~\ref{prop:polar} clarifies what is needed when $\beta>0$: if we show that $u_t<0$ and $0<\Ss<\Rr\tan\varphi$ in $\D$, then $P_1>0$ and $P_2>0$ immediately; we must then bound $P_1$ quantitatively from below in terms of $\ub^2$ to dominate the negative portion of $B_\beta$.

\section{Monotonicity and Hessian signs of the ground state}
\label{sec:directions}

Throughout Sections~\ref{sec:directions}--\ref{sec:close}, we restrict to $0<\theta\le\pi/4$ (so $0<q=\tan\theta\le1$ and $\cos(2\theta)\ge0$) and assume $\beta>0$, whence $\mu>0$. On the sloping side $\Gamma_\theta$, we orient the unit tangent and outward normal as
\[
 \tau=(\cos\theta,\sin\theta),\qquad n=(-\sin\theta,\cos\theta).
\]
By Lemma~\ref{lem:corner-regularity}, positive and negative parts of first and second partial derivatives of $u$ belong to $H^1(\D)$ and may be inserted directly into the ground-state identity \eqref{eq:ground-transform}.

\begin{lemma}[Cartesian monotonicity]\label{lem:first-directions}
In $\D$, one has $u_t<0$ and $u_s<0$. Moreover, $u_t<0$ on the open edges $\Gamma_R$ and $\Gamma_\theta$.
\end{lemma}

\begin{proof}
Let $p=u_t$. Differentiating \eqref{eq:unit-pde} gives $(-\Delta-\mu)p=0$ in $\D$, with $p_s+\beta p=0$ on $\Gamma_R$ and $p=0$ on $\Gamma_0$. Along $\Gamma_\theta$, $u_n=0$ implies $p=\sin\theta\,u_\tau$ and $u_{\tau n}=0$; differentiating $u_n=0$ tangentially and using $-\Delta u=\mu u$ yields
\begin{equation}\label{eq:p-sloping-boundary}
 p_n = \sin\theta\,u_{\tau n}+\cos\theta\,u_{nn}
     = -\cos\theta\,(u_{\tau\tau}+\mu u)
     = -\cot\theta\,p_\tau-\mu\cos\theta\,u \qquad\text{on }\Gamma_\theta.
\end{equation}
By \eqref{eq:corner-values}, $p(0,0)=0$ and $p(1,q)=-\beta q\,u(1,q)<0$, so $p_+=\max(p,0)$ vanishes at both endpoints of $\Gamma_\theta$. Testing $(-\Delta-\mu)p=0$ against $p_+$ and integrating $p_+ p_\tau = \frac12(p_+^2)_\tau$ along $\Gamma_\theta$ gives
\[
 \Q[p_+] = -\cot\theta\int_{\Gamma_\theta}p_+ p_\tau\dd\ell - \mu\cos\theta\int_{\Gamma_\theta}u p_+\dd\ell \le 0.
\]
By \eqref{eq:ground-transform}, $\Q[p_+]=0$, so $p_+$ is a constant multiple of $u$. Since $p_+=0$ on $\Gamma_0$, we must have $p_+\equiv0$, i.e., $p\le0$ in $\D$. If $p\equiv0$, then $u_\tau=0$ on $\Gamma_\theta$ together with $u_n=0$ would force $u_s=0$ on $\Gamma_\theta$ and hence $u$ would be constant in $t$ and satisfy $u_s(s,0)=0$, contradicting $u_s=-\beta u<0$ on $\Gamma_R$. Thus $p<0$ in $\D$ by the strong maximum principle. On $\Gamma_R$, a zero of $p$ would satisfy $p_s=-\beta p=0$, violating the Hopf lemma; on $\Gamma_\theta$, a zero of $p$ would be a boundary maximum with $p_\tau=0$, where \eqref{eq:p-sloping-boundary} gives $p_n=-\mu\cos\theta\,u<0$, again contradicting the Hopf lemma.

Next put $r=u_s$, which solves $(-\Delta-\mu)r=0$ in $\D$ with $\partial_n r=-r_t=0$ on $\Gamma_0$. On $\Gamma_R$, $r=-\beta u<0$; on $\Gamma_\theta$, $u_n=0$ gives $r=\cot\theta\,p\le0$. Hence $r_+$ vanishes on $\Gamma_R\cup\Gamma_\theta$, and testing with $r_+$ yields $\Q[r_+]=0$. Since $r_+=0$ on $\Gamma_R$, \eqref{eq:ground-transform} forces $r_+\equiv0$, and the strong maximum principle with $r|_{\Gamma_R}<0$ gives $u_s<0$ in $\D$.
\end{proof}

\begin{lemma}[Angular monotonicity]\label{lem:angular}
The rotational derivative $W:=s u_t - t u_s$ satisfies $W>0$ in $\D$. In addition,
\begin{equation}\label{eq:gaussian-moment-bound}
 \mu > 2\beta - \beta^2\left(\frac12+\frac{q^2}{6}\right).
\end{equation}
\end{lemma}

\begin{proof}
Set $\phi(s,t)=\exp[-\beta(s^2+t^2)/2]$. One checks immediately that $\phi$ satisfies $\phi_s+\beta\phi=0$ on $\Gamma_R$, $\partial_n\phi=0$ on $\Gamma_0\cup\Gamma_\theta$, and $-\Delta\phi=[2\beta-\beta^2(s^2+t^2)]\phi$. Applying Green's identity to $u$ and $\phi$ yields
\[
 \mu = 2\beta - \beta^2\frac{\int_\D(s^2+t^2)\phi u\dd s\dd t}{\int_\D\phi u\dd s\dd t}.
\]
By Lemma~\ref{lem:first-directions}, the product $w=\phi u$ satisfies $w_s<0$ and $w_t<0$ in $\D$. For each $s\in(0,1)$, Chebyshev's integral inequality comparing the strictly increasing function $t^2$ and the strictly decreasing function $t\mapsto w(s,t)$ on $(0,qs)$ gives
\[
 \int_0^{qs}t^2 w(s,t)\dd t < \frac{q^2 s^2}{3}\int_0^{qs}w(s,t)\dd t.
\]
Likewise, $m(s):=\int_0^1 w(s,qsz)\dd z$ is strictly decreasing on $(0,1)$, so Chebyshev's inequality with respect to the measure $s\dd s$ gives $\int_0^1 s^3 m(s)\dd s < \frac12\int_0^1 s m(s)\dd s$. Substituting $\int_\D w = q\int_0^1 s m(s)\dd s$ and $\int_\D s^2 w = q\int_0^1 s^3 m(s)\dd s$ proves \eqref{eq:gaussian-moment-bound}. Since $q\le1$, \eqref{eq:gaussian-moment-bound} implies in particular
\begin{equation}\label{eq:positive-angular-forcing}
 \mu+\beta^2-2\beta > \beta^2\left(\frac12-\frac{q^2}{6}\right) > 0.
\end{equation}

Because $s\partial_t-t\partial_s$ commutes with $\Delta$, we have $(-\Delta-\mu)W=0$ in $\D$, while $W=U_\varphi=r\Ss=0$ on both radial sides $\Gamma_0$ and $\Gamma_\theta$. Along $\Gamma_R$, differentiating $u_s=-\beta u$ in $t$ and using $u_{ss}=-\mu u-u_{tt}$ gives
\begin{equation}\label{eq:W-robin-boundary}
 W_s+\beta W = (\mu+\beta^2)t u + (t u_t)_t \qquad\text{on }\Gamma_R.
\end{equation}
Let $z=W_-=\max(-W,0)$, which vanishes on $\Gamma_0\cup\Gamma_\theta$ and hence at both endpoints $t=0$ and $t=q$ of $\Gamma_R$. On $\{t\in(0,q):z(1,t)>0\}$, we have $W(1,t)=-z(1,t)$, so $u_t(1,t)=W(1,t)+t u_s(1,t)=-z(1,t)-\beta t u(1,t)$. Testing $(-\Delta-\mu)W=0$ against $-z$ gives $\Q[z]=-\int_0^q [(\mu+\beta^2)t u+(t u_t)_t]z\dd t$; integrating $\int_0^q (t u_t)_t z\dd t = -\int_0^q t u_t z_t\dd t$ by parts, substituting $u_t=-z-\beta t u$ on $\{z>0\}$, and integrating $\int_0^q t z z_t\dd t$ and $\int_0^q t^2 u z_t\dd t$ by parts once more yields
\begin{equation}\label{eq:W-negative-part}
 \Q[z] = \int_0^q\left(\frac12-\beta t^2\right)z(1,t)^2\dd t
        -\int_0^q t\bigl(\mu+\beta^2-2\beta+\beta^2 t^2\bigr)u(1,t)z(1,t)\dd t
       \le \frac12\int_0^q z(1,t)^2\dd t,
\end{equation}
where the inequality follows from \eqref{eq:positive-angular-forcing} and $u,z\ge0$. On the other hand, $v:=z/u\in H^1(\D)$ vanishes on $\Gamma_\theta$, i.e., at $s=t/q$ for every $t\in(0,q)$. Since $u_s<0$ implies $u(s,t)\ge u(1,t)$ for $s\in[t/q,1]$, Cauchy--Schwarz gives
\[
 z(1,t)^2 = u(1,t)^2\left(\int_{t/q}^1 v_s(s,t)\dd s\right)^2
          \le u(1,t)^2\left(\int_{t/q}^1 u^2 v_s^2\dd s\right)\left(\int_{t/q}^1 u^{-2}\dd s\right)
          \le \int_{t/q}^1 u(s,t)^2 v_s(s,t)^2\dd s.
\]
Integrating in $t\in(0,q)$ and applying \eqref{eq:ground-transform} gives $\Q[z]\ge\int_0^q z(1,t)^2\dd t$, which combined with \eqref{eq:W-negative-part} forces $z(1,t)\equiv0$. Then $\Q[z]\le0$, so $z\equiv0$ and $W\ge0$ in $\D$. If $W\equiv0$, then $u_t=-\beta t u$ on $\Gamma_R$, and substituting this into \eqref{eq:W-robin-boundary} gives $0=(\mu+\beta^2-2\beta-\beta^2 t^2)t u$ on $(0,q)$, contradicting \eqref{eq:positive-angular-forcing} as $t\downarrow0$. Hence $W>0$ in $\D$ by the strong maximum principle.
\end{proof}

In terms of the polar gradient components $\Rr=-U_r$ and $\Ss=U_\varphi/r$ from Section~\ref{sec:polar}, Lemmas~\ref{lem:first-directions} and~\ref{lem:angular} translate immediately into
\begin{equation}\label{eq:RS-signs}
 \Rr>0 \quad\text{and}\quad 0<\Ss<\Rr\tan\varphi \qquad\text{in }\D,
\end{equation}
since $r\Ss=W>0$ and $u_t=-\Rr\sin\varphi+\Ss\cos\varphi<0$ (which also forces $\Rr>0$). In particular, both bulk integrals $P_1$ and $P_2$ in Proposition~\ref{prop:polar} are strictly positive for every $\theta\in(0,\pi/4]$ and $\beta>0$.

\begin{lemma}[Second-order derivative signs]\label{lem:hessian}
For $0<\theta\le\pi/4$ and $\beta>0$, one has in $\D$:
\[
 u_{st}>0,\qquad Z:=u_{ss}-u_{tt}\ge0,\qquad -u_{tt}\ge\frac{\mu}{2}u.
\]
More precisely, $Z>0$ in $\D$ when $\theta<\pi/4$, and $Z\equiv0$ when $\theta=\pi/4$.
\end{lemma}

\begin{proof}
Put $v=u_{st}$, which solves $(-\Delta-\mu)v=0$ in $\D$. On $\Gamma_R$, $v=-\beta u_t>0$ by Lemma~\ref{lem:first-directions}; on $\Gamma_0$, $v=0$. Along $\Gamma_\theta$, expressing $\partial_s,\partial_t$ in the $(\tau,n)$ frame using $u_n=u_{\tau n}=u_{\tau\tau n}=0$ and $u_{nn}=-\mu u-u_{\tau\tau}$ gives
\[
 v=\sin\theta\cos\theta\,(2u_{\tau\tau}+\mu u),\qquad
 v_n=-\cot(2\theta)\,v_\tau-\frac{\mu\cos(2\theta)}{2}u_\tau \qquad\text{on }\Gamma_\theta.
\]
By \eqref{eq:corner-values}, $v(0,0)=0$ and $v(1,q)=\beta^2 q\,u(1,q)>0$, so $v_-=\max(-v,0)$ vanishes near $(1,q)$ and at $(0,0)$. Testing $(-\Delta-\mu)v=0$ with $-v_-$ yields
\[
 \Q[v_-]=\frac{\mu\cos(2\theta)}{2}\int_{\Gamma_\theta}u_\tau v_-\dd\ell\le0,
\]
because $u_\tau=\csc\theta\,u_t<0$ on $\Gamma_\theta$ and $\cos(2\theta)\ge0$ for $\theta\le\pi/4$. Thus $v_-\equiv0$ by \eqref{eq:ground-transform} (as $v_-=0$ on $\Gamma_0$), and $v|_{\Gamma_R}>0$ implies $u_{st}>0$ in $\D$.

Next, $Z=u_{ss}-u_{tt}=-\mu u-2u_{tt}$ solves $(-\Delta-\mu)Z=0$ in $\D$, satisfies $Z_s+\beta Z=0$ on $\Gamma_R$ and $\partial_n Z=0$ on $\Gamma_0$, and has trace $Z=\cos(2\theta)(u_{\tau\tau}-u_{nn})=2\cot(2\theta)\,v\ge0$ on $\Gamma_\theta$. Hence $Z_-$ vanishes on $\Gamma_\theta$ and satisfies $\Q[Z_-]=0$, which forces $Z_-\equiv0$ by \eqref{eq:ground-transform}. If $\theta<\pi/4$, then $Z|_{\Gamma_\theta}>0$ gives $Z>0$ in $\D$; if $\theta=\pi/4$, then $Z|_{\Gamma_\theta}=0$ gives $\Q[Z]=0$ and hence $Z\equiv0$. Finally, $-u_{tt}=(\mu u+Z)/2\ge(\mu/2)u$.
\end{proof}

\section{Transverse power concavity and the boundary reference weight}
\label{sec:concavity}

To compare the Robin boundary trace $\ub(t)=u(1,t)$ against an explicit profile, we establish concavity of a suitable power of $u$ in the vertical variable $t$. By Lemma~\ref{lem:corner-regularity}, $u\in C^2(\overline\D)$ is smooth up to each open edge and strictly positive on $\overline\D$, so logarithmic derivatives of $u$ are $C^2$ in $\D\cup\Gamma_R\cup\Gamma_0\cup\Gamma_\theta$ and continuous on $\overline\D$.

\begin{proposition}[Transverse power concavity]\label{prop:power}
For $0<\theta\le\pi/4$ and $\beta>0$, set
\[
 a:=2\cos^2\theta=\frac{2}{1+q^2}\in[1,2).
\]
Then $\partial_{tt}(u^a)\le0$ on $\overline\D$. In particular, $t\mapsto\ub(t)^a$ is concave on $[0,q]$.
\end{proposition}

\begin{proof}
When $\theta=\pi/4$, we have $a=1$, and Lemma~\ref{lem:hessian} gives $\partial_{tt}u = -\frac12\mu u < 0$ directly. Assume henceforth that $0<\theta<\pi/4$, so $1<a<2$ and $\chi:=\cos(2\theta)>0$.

Write $v=\log u$, $H=v_{tt}+a v_t^2$, and introduce the elliptic operator $L=\Delta+2\nabla v\cdot\nabla$. Since $\partial_{tt}(u^a)=a u^a H$, it is equivalent to prove $H\le0$ on $\overline\D$. From $-\Delta u=\mu u$ we have $\Delta v+|\nabla v|^2=-\mu$; differentiating in $t$ gives $L v_t=0$ and $L v_{tt}=-2|\nabla v_t|^2$, whence
\begin{equation}\label{eq:power-subsolution}
 L H = 2(a-1)|\nabla v_t|^2 \ge 0 \qquad\text{in }\D.
\end{equation}

We examine $H$ on the three edges and three vertices of $\D$:
\begin{itemize}
\item On $\Gamma_R$, $v_s=-\beta$ is constant, so $v_{tts}=0$ and $v_{ts}=0$, giving
\begin{equation}\label{eq:power-robin-edge}
 H_s = 0 \qquad\text{on }\Gamma_R.
\end{equation}
\item On $\Gamma_0$ (including $(1,0)$ by continuity), $v_t=0$ and Lemma~\ref{lem:hessian} give
\begin{equation}\label{eq:power-axis}
 H = v_{tt} = \frac{u_{tt}}{u} \le -\frac{\mu}{2} < 0 \qquad\text{on }\Gamma_0\cup\{(1,0)\}.
\end{equation}
\item At the vertices $(0,0)$ and $(1,q)$, \eqref{eq:corner-values} gives $\nabla v(0,0)=0$ and $v_{tt}(0,0)=-\mu/2$, so $H(0,0)=-\mu/2<0$. At $(1,q)$, \eqref{eq:corner-values} yields $v_t=-\beta q$ and $u_{st}=\beta^2 q u$; differentiating $u_n=0$ tangentially along $\Gamma_\theta$ gives $(1-q^2)u_{st}=q(u_{ss}-u_{tt})=qZ$, so $Z(1,q)=\beta^2(1-q^2)u(1,q)$. Substituting $u_{tt}=-(\mu u+Z)/2$ and $a=2/(1+q^2)$ produces
\begin{equation}\label{eq:power-corner}
 H(1,q) = -\frac{\mu}{2}-\frac{\beta^2(1-q^2)}{2}+\frac{2\beta^2 q^2}{1+q^2} - \beta^2 q^2
        = -\frac{\mu}{2}-\frac{\beta^2(1-q^2)^2}{2(1+q^2)} < 0,
\end{equation}
where we note that $v_{tt}+a v_t^2 = \frac{u_{tt}}{u}+(a-1)\frac{u_t^2}{u^2}$.
\item Along the open sloping edge $\Gamma_\theta$, write $S_\theta=\sin\theta$ and $C_\theta=\cos\theta$. In the orthonormal coordinates $(\tau,n)$, the Neumann condition $v_n=0$ implies $v_{\tau n}=v_{\tau\tau n}=0$, and differentiating $\Delta v+|\nabla v|^2=-\mu$ in the normal direction gives $v_{nnn}=-v_{\tau\tau n}-2v_\tau v_{\tau n}-2v_n v_{nn}=0$. Expressing $\partial_t=S_\theta\partial_\tau+C_\theta\partial_n$ and eliminating $v_{nn}=-\mu-v_\tau^2-v_{\tau\tau}$, we find on $\Gamma_\theta$:
\[
 H = S_\theta^2 v_{\tau\tau}+C_\theta^2 v_{nn}+a S_\theta^2 v_\tau^2
   = -\chi v_{\tau\tau}-\mu C_\theta^2+(a S_\theta^2-C_\theta^2)v_\tau^2.
\]
Differentiating $H$ tangentially and normally along $\Gamma_\theta$ using $v_n=v_{\tau n}=v_{\tau\tau n}=v_{nnn}=0$ gives, for general $a$,
\begin{align*}
 \partial_\tau H &= -\chi v_{\tau\tau\tau} + 2(a S_\theta^2-C_\theta^2)v_\tau v_{\tau\tau},\\
 \partial_n H &= 2S_\theta C_\theta v_{\tau nn} + C_\theta^2 v_{nnn} + 2a S_\theta C_\theta v_\tau v_{nn}
               = -2S_\theta C_\theta\bigl[v_{\tau\tau\tau}+(a+2)v_\tau v_{\tau\tau}+a\mu v_\tau+a v_\tau^3\bigr],
\end{align*}
since $v_{\tau nn}=-\partial_\tau(\mu+v_\tau^2+v_{\tau\tau})=-2v_\tau v_{\tau\tau}-v_{\tau\tau\tau}$. Forming $\partial_n H-\tan(2\theta)\partial_\tau H$ cancels $v_{\tau\tau\tau}$; when $a=2C_\theta^2$, the cubic term $v_\tau^3$ also cancels against $v_\tau H$ because $a S_\theta^2-C_\theta^2 = C_\theta^2(2S_\theta^2-1) = -C_\theta^2\chi$, yielding the exact oblique relation
\begin{equation}\label{eq:power-sloping-edge}
 \partial_n H - \tan(2\theta)\,\partial_\tau H
 = \frac{4S_\theta C_\theta v_\tau}{\chi}\bigl(H+2\mu S_\theta^2 C_\theta^2\bigr)
 \qquad\text{on }\Gamma_\theta.
\end{equation}
\end{itemize}

Suppose for contradiction that $\max_{\overline\D}H>0$, attained at some $x_0\in\overline\D$. By \eqref{eq:power-subsolution} and the strong maximum principle, $x_0\notin\D$ (as $H$ cannot be a positive constant in view of \eqref{eq:power-axis}). By \eqref{eq:power-axis} and \eqref{eq:power-corner}, $x_0$ cannot lie on $\Gamma_0$ or at any of the three vertices. If $x_0\in\Gamma_R$, the Hopf lemma applied to \eqref{eq:power-subsolution} requires $H_s(x_0)>0$, contradicting \eqref{eq:power-robin-edge}. Finally, if $x_0\in\Gamma_\theta$, then $\partial_\tau H(x_0)=0$ and $H(x_0)>0$; since $v_\tau = u_\tau/u = \csc\theta\,(u_t/u) < 0$ on $\Gamma_\theta$ by Lemma~\ref{lem:first-directions}, \eqref{eq:power-sloping-edge} gives $\partial_n H(x_0)<0$, again contradicting the Hopf lemma. Therefore $H\le0$ on $\overline\D$.
\end{proof}

\begin{corollary}[Monotonicity of the boundary density ratio]\label{cor:weight}
Put $p:=1+q^2=2/a\in(1,2]$ and define
\begin{equation}\label{eq:w-def}
 w(t):=\left[1+\frac{2\beta q(q-t)}{p}\right]^p, \qquad 0\le t\le q.
\end{equation}
Then $t\mapsto\ub(t)^2/w(t)$ is positive and nondecreasing on $[0,q]$.
\end{corollary}

\begin{proof}
Let $V(t)=\ub(t)^a$ for $t\in[0,q]$. By \eqref{eq:corner-values}, $V'(q)=a\ub(q)^{a-1}u_t(1,q)=-a\beta q V(q)$, so the tangent line to $V$ at $t=q$ is
\[
 \ell(t)=V(q)\bigl[1+a\beta q(q-t)\bigr]>0.
\]
By Proposition~\ref{prop:power}, $V''\le0$ on $[0,q]$. Consequently $(V'\ell-V\ell')'=V''\ell\le0$, and since $(V'\ell-V\ell')(q)=0$, we have $V'\ell-V\ell'\ge0$ on $[0,q]$, meaning that $V/\ell$ is nondecreasing. Because $ap=2$ and $\ell(t)^p=V(q)^p w(t)$, $\ub(t)^2/w(t)=V(q)^{-p}(V(t)/\ell(t))^p$ is likewise nondecreasing.
\end{proof}

\section{Lower bound for the bulk term and boundary absorption (\texorpdfstring{$N\ge4$}{N >= 4})}
\label{sec:close}

We continue to assume $0<\theta\le\pi/4$ ($0<q\le1$) and $\beta>0$. Set $d:=1-c>0$, and define on $[0,q]$ the scalar kernels
\begin{equation}\label{eq:scalar-kernels}
 G(t)=t-c(1+t^2)\arctan t,\qquad
 K(t)=2ct\arctan t-\frac{d}{2},\qquad
 f(t)=tG(t).
\end{equation}
Note that $G(0)=G(q)=0$, $K(t)=\frac{d}{2}-G'(t)$, and $G(t)>0$ on $(0,q)$ (under $t=\tan\varphi$, $G(\tan\varphi)=\sec^2\varphi\,g(\varphi)$). In this notation, $B_\beta=\beta\int_0^q K(t)\ub(t)^2\dd t$.

We first bound the bulk term $P_1$ from below by an integral against $f(t)\ub(t)^2$. Put
\[
 k=\sqrt{\frac{\mu}{2}},\qquad
 C_k(s)=\cos\bigl(k(1-s)\bigr)+\frac{\beta}{k}\sin\bigl(k(1-s)\bigr),\qquad
 M_3=\int_0^1 s^3 C_k(s)^2\dd s.
\]

\begin{lemma}[One-dimensional bulk comparison]\label{lem:volume-lower}
For $0<\theta\le\pi/4$ and $\beta>0$, the bulk integral $P_1$ in \eqref{eq:P1} satisfies
\begin{equation}\label{eq:volume-lower}
 P_1 \ge \mu^2 M_3\int_0^q f(t)\ub(t)^2\dd t.
\end{equation}
\end{lemma}

\begin{proof}
Testing the Rayleigh quotient on $\D$ with $v\equiv1$ and with $v(s,t)=1-s$ gives $\mu\le 2\beta$ and $\mu\le 6$, respectively. Hence $k\le\sqrt3<\pi$ and $\beta\ge k^2$. Since $k\le\sqrt3$ implies $\tan(k/2)\le\frac{k/2}{1-k^2/8}\le\frac{4k}{5}<k$, we have for $z=k(1-s)\in[0,k]$:
\[
 C_k(s)-1 = \sin z\left(\frac{\beta}{k}-\tan\frac{z}{2}\right) \ge \sin z\,\bigl(k-\tan(k/2)\bigr) \ge 0.
\]

By Lemma~\ref{lem:hessian}, $u_{ss}+k^2 u = Z/2 \ge 0$ in $\D$. Fixing $t_0\in(0,q)$ and solving this ordinary differential inequality backward from $s=1$ (where $u(1,t_0)=\ub(t_0)$ and $u_s(1,t_0)=-\beta\ub(t_0)$) by variation of constants gives, for $t_0/q\le s\le1$,
\[
 u(s,t_0)=C_k(s)\ub(t_0)+\frac{1}{2k}\int_s^1\sin\bigl(k(r-s)\bigr)Z(r,t_0)\dd r \ge C_k(s)\ub(t_0),
\]
because $0\le k(r-s)\le k<\pi$ makes the sine kernel nonnegative. Evaluating at $t_0=st$ ($0\le s\le1$, $0\le t\le q$) and using $u_t<0$ together with $C_k(s)\ge1$ yields
\begin{equation}\label{eq:CG-ray-lower}
 u(s,st) \ge C_k(s)\ub(st) \ge C_k(s)\ub(t).
\end{equation}
Meanwhile, since $u_t(s,0)=0$, Lemma~\ref{lem:hessian} and $u_t<0$ give
\[
 -u_t(s,y) = \int_0^y (-u_{tt})(s,r)\dd r \ge \frac{\mu}{2}\int_0^y u(s,r)\dd r \ge \frac{\mu y}{2}u(s,y).
\]
Changing variables via $t=\tan\varphi$ in \eqref{eq:P1} rewrites $P_1$ as
\[
 P_1 = 2\mu\int_0^q G(t)\int_0^1 s^2 u(s,st)\bigl(-u_t(s,st)\bigr)\dd s\dd t.
\]
Inserting $-u_t(s,st)\ge\frac12\mu st\,u(s,st)$ and \eqref{eq:CG-ray-lower} yields \eqref{eq:volume-lower}.
\end{proof}

\begin{lemma}[Uniform moment bound]\label{lem:kernel-bound}
For every $\theta\in(0,\pi/4]$ and every $\beta>0$,
\begin{equation}\label{eq:kernel-bound}
 \mu^2 M_3 \ge \frac{\pi^4}{384}\beta^2 > \frac{\beta^2}{4}.
\end{equation}
\end{lemma}

\begin{proof}
Let $\kappa\in(0,\pi/2)$ be the unique root of $\kappa\tan\kappa=\beta$, and set $\phi(s,t)=\cos(\kappa s)\cos(\kappa t)$. On $\D$, $\phi>0$, $-\Delta\phi=2\kappa^2\phi$, $\phi_s+\beta\phi=0$ on $\Gamma_R$, and $\phi_t=0$ on $\Gamma_0$. Along $\Gamma_\theta$ ($t=qs$ with $0<q\le1$), convexity of $\tan x$ on $[0,\pi/2)$ implies $\tan(\kappa qs)\le q\tan(\kappa s)$, so
\[
 \partial_n\phi = \kappa\phi\bigl[\sin\theta\tan(\kappa s)-\cos\theta\tan(\kappa qs)\bigr] \ge 0 \qquad\text{on }\Gamma_\theta.
\]
Green's identity $(\mu-2\kappa^2)\int_\D u\phi = \int_{\Gamma_\theta}u\partial_n\phi\ge0$ then gives $\mu\ge 2\kappa^2$, i.e., $k\ge\kappa$.

\emph{Case 1: $k\le\pi/2$.} Then $\beta=\kappa\tan\kappa\le k\tan k$ (interpreted by continuity at $k=\pi/2$), which implies
\[
 \frac{C_k(s)}{\beta} = \frac{\cos(k(1-s))}{\beta}+\frac{\sin(k(1-s))}{k}
 \ge \frac{\cos(k(1-s))}{k\tan k}+\frac{\sin(k(1-s))}{k} = \frac{\cos(ks)}{k\sin k}.
\]
Since $\mu=2k^2$, we obtain
\[
 \frac{\mu^2 M_3}{\beta^2} \ge F(k) := \frac{4k^2}{\sin^2 k}\int_0^1 s^3\cos^2(ks)\dd s.
\]
Because $\frac{\dd^2}{\dd z^2}\cos(k\sqrt z)=\frac{k}{4z^{3/2}}(\sin(k\sqrt z)-k\sqrt z\cos(k\sqrt z))\ge0$ for $z\in(0,1]$, the function $z\mapsto\cos(k\sqrt z)$ is convex on $[0,1]$, so its tangent at $z=1$ gives $\cos(ks)\ge\cos k+\frac12 k\sin k\,(1-s^2)\ge0$. Squaring and integrating against $s^3\dd s$ yields, with $y=k\cot k$,
\[
 F(k) \ge y^2+\frac{k^2 y}{3}+\frac{k^4}{24}.
\]
With $\gamma=4/\pi^2$, we claim that $y\ge 1-\gamma k^2$ on $[0,\pi/2]$. Indeed, $\psi(k):=k\cos k-(1-\gamma k^2)\sin k$ vanishes at $k=0$ and $k=\pi/2$, and $\psi'(k)=k\sin k[(2\gamma-1)+\gamma k\cot k]$ changes sign once from positive to negative on $(0,\pi/2)$, so $\psi\ge0$. Substituting $y\ge 1-\gamma k^2$ gives $F(k)\ge P(k^2)$, where
\[
 P(x)=1+\left(\frac13-2\gamma\right)x+\left(\gamma^2-\frac\gamma3+\frac1{24}\right)x^2.
\]
Here $\gamma^2-\frac\gamma3+\frac1{24}=(\gamma-\frac16)^2+\frac1{72}>0$ while $P'(1/\gamma)=-\frac13+\frac{1}{12\gamma}=\frac{\pi^2-16}{48}<0$, so $P$ decreases on $[0,\pi^2/4]$ and $F(k)\ge P(\pi^2/4)=\pi^4/384$.

\emph{Case 2: $k>\pi/2$.} Setting $x=1-s$ and restricting the integral for $M_3$ to $0\le x\le\pi/(2k)$, where $\cos(kx),\sin(kx)\ge0$ and hence $C_k(1-x)\ge\frac{\beta}{k}\sin(kx)$, we find
\[
 \frac{\mu^2 M_3}{\beta^2} \ge 4k^2\int_0^{\pi/(2k)}(1-x)^3\sin^2(kx)\dd x
 = 4\int_0^{\pi/2}\frac{(k-z)^3}{k^2}\sin^2 z\dd z.
\]
For each $z\in[0,\pi/2]$, $\frac{\dd}{\dd k}\bigl[(k-z)^3/k^2\bigr]=(k-z)^2(k+2z)/k^3\ge0$ for $k\ge\pi/2$. Evaluating at $k=\pi/2$ recovers $F(\pi/2)\ge\pi^4/384$. Since $\pi^4>97.4>96$, \eqref{eq:kernel-bound} follows.
\end{proof}

We now show that the lower bound $\frac{\beta^2}{4}\int_0^q f\ub^2$ absorbs $B_\beta$. The transfer from $\ub^2$ to the explicit weight $w$ of Corollary~\ref{cor:weight} relies on three elementary scalar lemmas.

\begin{lemma}[Single-crossing transfer]\label{lem:one-crossing}
For every $L\ge0$ and $0<q\le1$, the function $K+Lf$ has a unique zero $t_0\in(0,q)$ and changes sign at $t_0$ from negative to positive. More generally, if $h$ is continuous on $[a_0,b_0]$ with $h\le0$ on $[a_0,t_0]$ and $h\ge0$ on $[t_0,b_0]$, while $w>0$ is integrable and $R>0$ is nondecreasing, then
\begin{equation}\label{eq:one-crossing-transfer}
 \int_{a_0}^{b_0} h(t)R(t)w(t)\dd t \ge R(t_0)\int_{a_0}^{b_0} h(t)w(t)\dd t.
\end{equation}
\end{lemma}

\begin{proof}
Clearly $K$ is strictly increasing on $[0,q]$, with $K(0)=-d/2<0$ and $K(q)=2\sin^2\theta-d/2>0$ (since $c=\frac{\sin(2\theta)}{2\theta}\ge\cos(2\theta)=1-2\sin^2\theta$ implies $d\le 2\sin^2\theta$). Also $f\ge0$ on $[0,q]$. Whenever $t\in(0,q]$ satisfies $K(t)\le0$, we have $ct\arctan t\le d/4$, and hence
\[
 f'(t)=(2-c)t-c\arctan t-3ct^2\arctan t \ge 2dt-3ct^2\arctan t \ge \frac54 dt > 0,
\]
using $\arctan t\le t$. Thus $K+Lf$ is strictly increasing on $\{K\le0\}$ and strictly positive on $\{K>0\}$ (when $L>0$; for $L=0$, $K$ itself is strictly increasing), proving the unique sign change from negative to positive. Finally, \eqref{eq:one-crossing-transfer} follows by integrating the pointwise inequality $h(t)[R(t)-R(t_0)]w(t)\ge0$.
\end{proof}

\begin{lemma}[Two-sided polynomial bounds for $G$]\label{lem:G-bounds}
For $0\le t\le q\le1$,
\begin{equation}\label{eq:G-two-sided}
 \frac78 d\,t\left(1-\frac{t^2}{q^2}\right)
 \le G(t) \le
 d\,t\left(1-\frac{t^2}{q^2}\right).
\end{equation}
\end{lemma}

\begin{proof}
Define on $[0,1]$ the function
\[
 \Phi(z):=\frac{(1+z)\arctan\sqrt z}{\sqrt z}=\int_0^1\frac{1+z}{1+zy^2}\dd y,\qquad \Phi(0)=1,
\]
which satisfies $\Phi'(z)=\int_0^1\frac{1-y^2}{(1+zy^2)^2}\dd y>0$ and $\Phi''(z)<0$. Since $c=1/\Phi(q^2)$ and $d=1-c=c[\Phi(q^2)-1]$, we have $G(t)/t = c[\Phi(q^2)-\Phi(t^2)]$. Concavity of $\Phi$ on $[0,q^2]$ gives the secant upper bound $\Phi(q^2)-\Phi(t^2)\le[\Phi(q^2)-1](1-t^2/q^2)$ and the tangent lower bound $\Phi(q^2)-\Phi(t^2)\ge \Phi'(q^2)(q^2-t^2)$. This immediately yields the upper bound in \eqref{eq:G-two-sided}. For the lower bound, note that $(1-7x)/(1+x)^2$ decreases on $[0,1]$, so for every $z\in(0,1]$,
\[
 \Phi'(z)-\frac{7}{8z}\bigl(\Phi(z)-1\bigr)
 =\frac18\int_0^1(1-y^2)\frac{1-7zy^2}{(1+zy^2)^2}\dd y
 \ge \Phi'(1)-\frac78\bigl(\Phi(1)-1\bigr)
 =\frac{22-7\pi}{16}>0,
\]
since $\Phi(1)=\pi/2$, $\Phi'(1)=1/2$, and $\pi<22/7$. Evaluating at $z=q^2$ gives $q^2\Phi'(q^2)\ge\frac78[\Phi(q^2)-1]$, proving the lower bound in \eqref{eq:G-two-sided}.
\end{proof}

\begin{lemma}[Positivity against the reference weight]\label{lem:scalar}
For $0<q\le1$, $p=1+q^2$, $\beta\ge0$, and $w(t)$ given by \eqref{eq:w-def},
\begin{equation}\label{eq:reference-boundary}
 \int_0^q\left[K(t)+\frac{\beta}{4}f(t)\right]w(t)\dd t > 0.
\end{equation}
\end{lemma}

\begin{proof}
Set $B=2\beta q^2/p\ge0$, so $w(t)=[1+B(1-t/q)]^p$ and $w'\le0$. Using $K=d/2-G'$ and $G(0)=G(q)=0$, an integration by parts gives
\[
 \int_0^q\left[K(t)+\frac{\beta}{4}f(t)\right]w(t)\dd t
 =\frac{d}{2}\int_0^q w(t)\dd t + \int_0^q G(t)w'(t)\dd t + \frac{\beta}{4}\int_0^q t G(t)w(t)\dd t.
\]
Since $w'\le0$ and $w\ge0$, Lemma~\ref{lem:G-bounds} bounds this below by replacing $G(t)$ with $dt(1-t^2/q^2)$ in the middle integral and with $\frac78 dt(1-t^2/q^2)$ in the last integral. Integrating $\int_0^q t(1-t^2/q^2)w'(t)\dd t = -\int_0^q(1-3t^2/q^2)w(t)\dd t$ by parts and substituting $t=qx$ ($b:=pB=2\beta q^2\ge0$) yields the lower bound
\begin{equation}\label{eq:scalar-polynomial-reduction}
 dq\int_0^1 P_b(x)\left[1+\frac{b}{p}(1-x)\right]^p\dd x,
 \qquad P_b(x):=3x^2-\frac12+\frac{7b}{64}x^2(1-x^2).
\end{equation}
Viewed as a quadratic in $x^2$, $P_b$ is concave (or linear if $b=0$) with $P_b(0)=-1/2<0$ and $P_b(1)=5/2>0$, so $P_b$ changes sign once on $(0,1)$ from negative to positive. Put $w_*(x)=[1+\frac{b}{2}(1-x)]^2$. Because $1<p\le2$, the ratio $R_p(x):=[1+\frac{b}{p}(1-x)]^p/w_*(x)$ satisfies
\[
 \frac{\dd}{\dd x}\log R_p(x) = \frac{b}{1+\frac{b}{2}(1-x)} - \frac{b}{1+\frac{b}{p}(1-x)} \ge 0,
\]
so $R_p>0$ is nondecreasing on $[0,1]$. Applying \eqref{eq:one-crossing-transfer} reduces positivity of \eqref{eq:scalar-polynomial-reduction} to positivity of $\int_0^1 P_b(x)w_*(x)\dd x$. With $z=b/2\ge0$, direct integration gives
\[
 \int_0^1 P_b(x)w_*(x)\dd x = \frac{5z^3-43z^2+28z+480}{960}.
\]
For $0\le z\le5$, $5z^3-43z^2+480$ decreases from $480$ to $30$, while $28z\ge0$; for $z=5+y$ with $y\ge0$, the numerator equals $5y^3+32(y-\frac{27}{64})^2+\frac{21031}{128}>0$. Since $dq>0$, \eqref{eq:reference-boundary} holds.
\end{proof}

Combining these ingredients closes the positive Robin branch for $N\ge4$.

\begin{proposition}[Monotonicity for $\alpha>0$ and $N\ge4$]\label{prop:close-N4}
For every $A>0$ and $\alpha>0$, $\Lambda'(\theta)>0$ for all $\theta\in(0,\pi/4]$. Consequently, $\lambda_1^\alpha(P_{N+1,A})<\lambda_1^\alpha(P_{N,A})$ for every integer $N\ge4$.
\end{proposition}

\begin{proof}
Fix $\theta\in(0,\pi/4]$ and $\beta=\alpha h(\theta)>0$. By Corollary~\ref{cor:weight}, $R:=\ub^2/w$ is positive and nondecreasing on $[0,q]$, so Lemmas~\ref{lem:one-crossing} and~\ref{lem:scalar} give
\[
 \int_0^q\left[K(t)+\frac{\beta}{4}f(t)\right]\ub(t)^2\dd t > 0.
\]
Assembling Proposition~\ref{prop:polar}, \eqref{eq:RS-signs} (which gives $P_2>0$), Lemma~\ref{lem:volume-lower}, Lemma~\ref{lem:kernel-bound}, and $f\ge0$, we obtain
\[
 J = P_1+P_2+B_\beta
   \ge P_2 + \beta\int_0^q\left[K(t)+\frac{\mu^2 M_3}{\beta}f(t)\right]\ub(t)^2\dd t
   > \beta\int_0^q\left[K(t)+\frac{\beta}{4}f(t)\right]\ub(t)^2\dd t > 0.
\]
Equations~\eqref{eq:branch-derivative} and~\eqref{eq:J-unit} therefore give $\Lambda'(\theta)>0$ on $(0,\pi/4]$, and integrating over $[\pi/(N+1),\pi/N]$ proves $\lambda_1^\alpha(P_{N+1,A})<\lambda_1^\alpha(P_{N,A})$ for every $N\ge4$.
\end{proof}

\section{The equilateral triangle, the square, and the Dirichlet limit}
\label{sec:final}

To complete the proof of Theorem~\ref{thm:main}, it remains to establish $\lambda_1^\alpha(P_{3,A})>\lambda_1^\alpha(P_{4,A})$ for $\alpha>0$. Write
\[
 k_3=\sqrt{\lambda_1^\alpha(P_{3,A})}>0,\qquad
 k_4=\sqrt{\lambda_1^\alpha(P_{4,A})}>0,\qquad
 r=\frac{\sqrt A}{3^{3/4}},
\]
where $r$ is the inradius of $P_{3,A}$. Separation of variables on the square $[-\sqrt A/2,\sqrt A/2]^2$ shows that $k_4/\sqrt2$ is the smallest positive root of $\xi\tan(\xi\sqrt A/2)=\alpha$, which reads
\begin{equation}\label{eq:square-secular}
 \frac{k_4\sqrt A}{2\sqrt2} = \arctan\frac{\sqrt2\alpha}{k_4}.
\end{equation}

\begin{lemma}[Equilateral triangle versus square]\label{lem:triangle-square}
For every $A>0$ and $\alpha>0$, $k_3$ is the unique positive root of
\begin{equation}\label{eq:triangle-secular}
 \frac{3r k_3}{2} = \arctan\frac{\alpha}{k_3} + \arctan\frac{2\alpha}{k_3},
\end{equation}
and satisfies $\lambda_1^\alpha(P_{3,A})>\lambda_1^\alpha(P_{4,A})$.
\end{lemma}

\begin{proof}
Equation~\eqref{eq:triangle-secular} is the fundamental symmetric mode ($m=n=0$) of the secular system for the Robin equilateral triangle \cite[eqs.~(1.1)--(1.2)]{RudnickWigman}. To make the paper self-contained, we verify directly that its unique positive root $k$ gives $\lambda_1^\alpha(P_{3,A})=k^2$. Let $n_1,n_2,n_3\in\RR^2$ be unit vectors at $120^\circ$ angles, so that $P_{3,A}=\{z\in\RR^2:n_j\cdot z<r\text{ for }j=1,2,3\}$. For the unique $k>0$ solving $3rk/2=\arctan(\alpha/k)+\arctan(2\alpha/k)$, set
\[
 \hat a=\arctan\frac{\alpha}{k},\qquad \hat b=\arctan\frac{2\alpha}{k},\qquad \delta=\frac{\hat a-2\hat b}{3},
\]
and define $u_\alpha(z)=\sum_{j=1}^3\cos(k n_j\cdot z+\delta)$, which solves $-\Delta u_\alpha=k^2 u_\alpha$. Along any side $n_j\cdot z=r$, let $t\in[-\sqrt A/3^{1/4},\sqrt A/3^{1/4}]=[-\sqrt3 r,\sqrt3 r]$ be the centered arc-length coordinate; using $kr+\delta=\hat a$ and $-\frac12 kr+\delta=-\hat b$, we obtain on that side
\[
 u_\alpha = \cos\hat a + 2\cos\hat b\,\cos\frac{\sqrt3 kt}{2},\qquad
 \partial_\nu u_\alpha+\alpha u_\alpha
 = (-k\sin\hat a+\alpha\cos\hat a) + (-k\sin\hat b+2\alpha\cos\hat b)\cos\frac{\sqrt3 kt}{2} = 0.
\]
Thus $u_\alpha$ is a Robin eigenfunction with eigenvalue $k^2$. Moreover, along the side we have $|\sqrt3 kt/2|\le \frac32 rk = \hat a+\hat b < \pi$; since $\cos\hat b>0$, the boundary trace of $u_\alpha$ attains its minimum at the vertices:
\[
 \min_{\partial P_{3,A}}u_\alpha = \cos\hat a + 2\cos\hat b\,\cos(\hat a+\hat b) = \frac{3\cos\hat a}{1+4(\alpha/k)^2} > 0.
\]
To identify $u_\alpha$ with the principal eigenfunction on $P_{3,A}$, observe that as $\alpha\downarrow0$, $k(\alpha)\to0$ and $u_\alpha\to3>0$ uniformly on $\overline{P_{3,A}}$. Because the Robin eigenvalues vary continuously with $\alpha$ and the principal eigenvalue $\lambda_1^\alpha(P_{3,A})$ is simple for every $\alpha\ge0$ (with eigenfunctions of strict sign), the continuous branch $\alpha\mapsto k(\alpha)^2$ starts at $\lambda_1^0(P_{3,A})=0$ for $\alpha=0$ and can never cross a higher eigenvalue at any $\alpha>0$: the subset of $[0,\infty)$ where $k(\alpha)^2=\lambda_1^\alpha(P_{3,A})$ is nonempty, closed by continuity, and open by simplicity of $\lambda_1^\alpha$ and the local spectral gap. Hence $k_3=k(\alpha)$, proving \eqref{eq:triangle-secular}.

Now define, for $x>0$,
\[
 \Psi(x) := \arctan x + \arctan(2x) - \frac{4\sqrt2}{3}\arctan(\sqrt2 x).
\]
We have $\Psi(0)=0$ and
\[
 \Psi'(x) = \frac{1}{1+x^2}+\frac{2}{1+4x^2}-\frac{8}{3(1+2x^2)}
          = \frac{(2x^2-1)^2}{3(1+x^2)(1+2x^2)(1+4x^2)} \ge 0,
\]
with equality only at $x=1/\sqrt2$, so $\Psi(x)>0$ for all $x>0$. Since $(4/3)^4=256/81>3$, we also have $\frac{4\sqrt2}{3}>\sqrt2\,3^{1/4}$. Evaluating $\Psi$ at $x=\alpha/k_4>0$ and using \eqref{eq:square-secular} together with $r=\sqrt A/3^{3/4}$ gives
\[
 \arctan\frac{\alpha}{k_4}+\arctan\frac{2\alpha}{k_4}
 > \frac{4\sqrt2}{3}\arctan\frac{\sqrt2\alpha}{k_4}
 > \sqrt2\,3^{1/4}\arctan\frac{\sqrt2\alpha}{k_4}
 = \frac{3r k_4}{2}.
\]
Because $k\mapsto \arctan(\alpha/k)+\arctan(2\alpha/k)-\frac32 rk$ is strictly decreasing on $(0,\infty)$ and vanishes at $k_3$, it follows that $k_3>k_4$, i.e., $\lambda_1^\alpha(P_{3,A})>\lambda_1^\alpha(P_{4,A})$.
\end{proof}

\begin{proof}[Proof of Theorem~\ref{thm:main}]
Combine Proposition~\ref{prop:negative} (for $\alpha\le0$ and all $N\ge3$), Proposition~\ref{prop:close-N4} (for $\alpha>0$ and $N\ge4$), and Lemma~\ref{lem:triangle-square} (for $\alpha>0$ and $N=3$).
\end{proof}

\begin{proof}[Proof of Corollary~\ref{cor:dirichlet}]
For each $\theta\in(0,\pi/2)$, retain the triangle $D_h(\theta)$ of height $h(\theta)=(A\theta/(\pi\tan\theta))^{1/2}$ from Section~\ref{sec:branch}. Write $\Lambda_\alpha(\theta)$ for its first mixed Robin--Neumann eigenvalue ($\alpha>0$) and $\Lambda_D(\theta)$ for its first mixed Dirichlet--Neumann eigenvalue, with $u=0$ on $\Gamma_{R,h}=\{s=h(\theta)\}$ and $\partial_n u=0$ on the two radial sides. Reflection symmetry and simplicity give $\Lambda_D(\pi/N)=\lambda_1^D(P_{N,A})$ for every $N\ge3$.

For fixed $\theta\in(0,\pi/2)$, the Rayleigh quotient shows that $\alpha\mapsto\Lambda_\alpha(\theta)$ is nondecreasing. Restricting the trial functions to $V_\theta=\{v\in H^1(D_h(\theta)):\operatorname{Tr}_{\Gamma_{R,h}}v=0\}$ gives $\Lambda_\alpha(\theta)\le\Lambda_D(\theta)$. Taking $L^2$-normalized eigenfunctions $u_\alpha$ for $\Lambda_\alpha(\theta)$, the uniform bound $\int_{D_h}|\nabla u_\alpha|^2+\alpha\int_{\Gamma_{R,h}}u_\alpha^2\le\Lambda_D(\theta)$ yields a weakly $H^1$-convergent (and strongly $L^2$-convergent) subsequence as $\alpha\to+\infty$ whose limit has zero trace on $\Gamma_{R,h}$, so lower semicontinuity implies $\lim_{\alpha\to+\infty}\Lambda_\alpha(\theta)=\Lambda_D(\theta)$ (and similarly $\lim_{\alpha\to+\infty}\lambda_1^\alpha(P_{N,A})=\lambda_1^D(P_{N,A})$). By Proposition~\ref{prop:close-N4}, $\Lambda_\alpha$ is increasing on $(0,\pi/4]$ for each $\alpha>0$; passing to the limit shows that $\Lambda_D$ is nondecreasing on $(0,\pi/4]$.

To prove strictness, pull $D_h(\theta)$ back to $D_*=\{0<\eta<\xi<1\}$ via $(s,t)=(h(\theta)\xi,h(\theta)\tan\theta\,\eta)$. On the fixed form domain $V=\{v\in H^1(D_*):v(1,\cdot)=0\}$,
\[
 \Lambda_D(\theta) = \min_{v\in V\setminus\{0\}}
 \frac{\displaystyle\int_{D_*}\Bigl[h(\theta)^{-2}|v_\xi|^2+\bigl(h(\theta)\tan\theta\bigr)^{-2}|v_\eta|^2\Bigr]}{\displaystyle\int_{D_*}v^2}.
\]
The form coefficients $h(\theta)^{-2}$ and $(h(\theta)\tan\theta)^{-2}$ are real analytic and locally uniformly positive on $(0,\pi/2)$, and $\Lambda_D(\theta)$ is simple for every $\theta\in(0,\pi/2)$. By analytic perturbation theory (or the real-analytic implicit function theorem applied to the normalized weak equation on $V$, whose linearization is invertible on $u^\perp$ by the spectral gap), $\theta\mapsto\Lambda_D(\theta)$ is real analytic on $(0,\pi/2)$. On the other hand, on almost every horizontal slice $t\in(0,h\tan\theta)$ of $D_h(\theta)$, the segment $(t\cot\theta,h)$ has length at most $h$ and vanishes at $s=h$, so the one-dimensional mixed Poincar\'e inequality integrated in $t$ gives
\[
 \Lambda_D(\theta) \ge \frac{\pi^2}{4h(\theta)^2} = \frac{\pi^3\tan\theta}{4A\theta} \longrightarrow +\infty \qquad\text{as }\theta\uparrow\frac\pi2.
\]
Thus $\Lambda_D$ is nonconstant on $(0,\pi/2)$. If $\Lambda_D(\theta_1)=\Lambda_D(\theta_2)$ for some $0<\theta_1<\theta_2\le\pi/4$, monotonicity on $(0,\pi/4]$ would force $\Lambda_D$ to be constant on $[\theta_1,\theta_2]$, and real analyticity would force $\Lambda_D$ to be constant throughout $(0,\pi/2)$, a contradiction. Hence $\Lambda_D$ is strictly increasing on $(0,\pi/4]$, proving $\lambda_1^D(P_{N+1,A})<\lambda_1^D(P_{N,A})$ for all $N\ge4$.

For $N=3$, letting $\alpha\to+\infty$ in \eqref{eq:triangle-secular} and \eqref{eq:square-secular} gives $\lambda_1^D(P_{3,A})=\frac{4\pi^2}{\sqrt3 A}>\frac{2\pi^2}{A}=\lambda_1^D(P_{4,A})$. Finally, $P_{N,A}$ contains the concentric disk of inradius $r_N=\sqrt{A/\pi}\sqrt{(\pi/N)/\tan(\pi/N)}$ and is contained in the concentric disk of circumradius $R_N=r_N\sec(\pi/N)$. Since $r_N,R_N\to\sqrt{A/\pi}$ as $N\to\infty$, Dirichlet domain inclusion $\lambda_1^D(B_{\pi R_N^2})\le\lambda_1^D(P_{N,A})\le\lambda_1^D(B_{\pi r_N^2})$ and scaling yield $\lim_{N\to\infty}\lambda_1^D(P_{N,A})=\lambda_1^D(B_A)$.
\end{proof}

\section{The fixed-perimeter normalization}
\label{sec:perimeter}

\begin{proof}[Proof of Corollary~\ref{cor:perimeter}]
Write
\[
 A_N=|Q_{N,L}|=\frac{L^2}{4N\tan(\pi/N)}.
\]
Since $x\mapsto\tan x/x$ is strictly increasing on $(0,\pi/2)$ (by strict convexity of $\tan x$ and $\tan 0=0$), the sequence $A_N$ is strictly increasing and converges to $L^2/(4\pi)$.

Suppose first that $\alpha>0$. For any bounded connected Lipschitz domain $\Omega$ and any $r>1$, the Rayleigh quotient on the dilated domain is
\[
 \lambda_1^\alpha(r\Omega)=
 \min_{v\in H^1(\Omega)\setminus\{0\}}
 \frac{r^{-2}\int_\Omega|\nabla v|^2+\alpha r^{-1}\int_{\partial\Omega}v^2}
      {\int_\Omega v^2}
 \le r^{-1}\lambda_1^\alpha(\Omega)<\lambda_1^\alpha(\Omega),
\]
where the final strict inequality uses $\lambda_1^\alpha(\Omega)>0$. Applying this to the dilation from $P_{N+1,A_N}$ to $P_{N+1,A_{N+1}}$ and then using Theorem~\ref{thm:main}, we obtain
\[
 \lambda_1^\alpha(Q_{N+1,L})
 =\lambda_1^\alpha(P_{N+1,A_{N+1}})
 <\lambda_1^\alpha(P_{N+1,A_N})
 <\lambda_1^\alpha(P_{N,A_N})
 =\lambda_1^\alpha(Q_{N,L}).
\]
For $\alpha=0$, the conclusion follows from the constant eigenfunction.

Now suppose that $\alpha<0$. For $0<\theta\le\pi/3$, use the mixed problem of Section~\ref{sec:branch} on the right triangle of height and vertical side length
\[
 h=h_L(\theta)=\frac{L\theta}{2\pi\tan\theta},\qquad
 a=a_L(\theta)=h_L(\theta)\tan\theta=\frac{L\theta}{2\pi},
\]
and denote its principal eigenvalue by $\Lambda_L(\theta)$. At $\theta=\pi/N$, the triangle is a fundamental symmetry sector of $Q_{N,L}$, so $\Lambda_L(\pi/N)=\lambda_1^\alpha(Q_{N,L})$. As in Section~\ref{sec:branch}, the pullback to $D_*$ and differentiation of the quadratic form give
\[
 \Lambda_L'(\theta)=\left(\frac{a'}a-\frac{h'}h\right)(Y_h-X_h)
 +\alpha\frac{a'}a B_h-\Lambda_L(\theta)\left(\frac{h'}h+\frac{a'}a\right),
\]
where $X_h,Y_h,B_h$ are defined using the $L^2$-normalized ground state on this triangle. In the present normalization,
\[
 \frac{h'}h=\frac1\theta-\frac1{\sin\theta\cos\theta},
 \qquad \frac{a'}a=\frac1\theta.
\]
With $c=\sin\theta\cos\theta/\theta$ and $\Lambda_L=X_h+Y_h+\alpha B_h$, the derivative reduces to
\begin{equation}\label{eq:perimeter-derivative}
 \Lambda_L'(\theta)=\frac{2}{\sin\theta\cos\theta}
 \left[(1-c)\left(Y_h+\frac\alpha2 B_h\right)-cX_h\right].
\end{equation}
At each fixed $\theta$, scale the triangle to unit height as in Section~\ref{sec:branch}, with $\beta=\alpha h<0$. The Rellich identity established in the proof of Proposition~\ref{prop:negative} depends only on the unit-height mixed problem, not on the normalization of the original polygon. Hence
\[
 h^2\left(Y_h+\frac\alpha2 B_h\right)
 =Y+\frac\beta2 B
 =-\frac{\sin\theta}{2}\int_{\Gamma_\theta}(1-s)
     \bigl(|\nabla u|^2-\mu u^2\bigr)\dd\ell<0.
\]
Since $0<c<1$ and $X_h\ge0$, \eqref{eq:perimeter-derivative} gives $\Lambda_L'(\theta)<0$ on $(0,\pi/3]$. Integrating between $\pi/(N+1)$ and $\pi/N$ proves the strict inequality for every $N\ge3$.

Finally, by Dirichlet homogeneity and Corollary~\ref{cor:dirichlet},
\[
 \lambda_1^D(Q_{N,L})=\frac{1}{A_N}\lambda_1^D(P_{N,1}).
\]
Both positive factors on the right decrease strictly with $N$. Since $A_N\to L^2/(4\pi)$ and $\lambda_1^D(P_{N,1})\to\lambda_1^D(B_1)$, scaling also gives the limit in \eqref{eq:perimeter-dirichlet}.
\end{proof}

\section*{Declaration of generative AI use}
Large language models (specifically GPT-5.5 and Gemini 4 Carbon) were used during the preparation of this manuscript for purposes including idea exploration, drafting and language editing, and related editorial assistance. The authors carefully reviewed and verified the final manuscript and take full responsibility for its content, mathematical correctness, and conclusions.

\end{document}